# A Multi-Stage Supply Chain Network Optimization Using Genetic Algorithms

Nelson Christopher Dzupire[1*], Yaw Nkansah-Gyekye[1]

[1]Nelson Mandela African Institution Of Science and Engineering,

School of Computational and Communication Science and Engineering,

P. O. Box 447, Arusha, Tanzania.

* E-mail of corresponding author: dzupiren@nm-aist.ac.tz

**Abstract**

In today's global business market place, individual firms no longer compete as independent entities with unique brand names but as integral part of supply chain links. Key to success of any business is satisfying customer's demands on time which may result in cost reductions and increase in service level. In supply chain networks decisions are made with uncertainty about product's demands, costs, prices, lead times, quality in a competitive and collaborative environment. If poor decisions are made, they may lead to excess inventories that are costly or to insufficient inventory that cannot meet customer's demands.

In this work we developed a bi-objective model that minimizes system wide costs of the supply chain and delays on delivery of products to distribution centers for a three echelon supply chain. Picking a set of Pareto front for multi-objective optimization problems require robust and efficient methods that can search an entire space. We used evolutionary algorithms to find the set of Pareto fronts which have proved to be effective in finding the entire set of Pareto fronts.

## 1. Introduction

A supply chain (SC) is a set of facilities, suppliers, customers, products and methods of controlling inventory, purchasing and distribution (Sabri and Beamon, 2000). It is aimed at providing customers with products they want in a timely way and efficiently and as profitable as possible. The main objective is to enhance the operational efficiency, profitability and competitive position of a firm and its supplier chain partners (Min and Zhou, 2002). Each independent entity of supply chain has inherent objective function to maximize in business transactions for profit maximization (Pinto, 2004). In a supply chain network (SCN) managers need to make strategic decisions that are viable for the business. The decisions range from what product to produce and their design, how much, when and from where to buy a product, how much, where and when to produce a product, etc (Veinott Jr, 2002). These strategic decisions are made with uncertainty about product demands, costs, prices, lead times and quality. Besides, the environment in which these decisions are made is competitive.

Recently there has been a growing interest in research in supply chain network optimization problems. This may be due to increasing competitiveness introduced by rapid globalization such that firms wants to reduce costs and maintain profit margins as observed by Altiparmak *et al.* (2006) or from a practical stand point of view, the rise may be from a number of changes in the manufacturing environment including the rising cost of manufacturing, the shrinking of manufacturing bases, shortened life cycle, the leveling of the playing field within the manufacturing industries and globalization of market economies as suggested by Beamon (1998) as well as attractive to cost ratios and building of long term relationships with trusted suppliers as stated by Williams and Gunal (2003).

Supply chain network design problems are functions of different parameters namely lead times at each entity, fill rates, inventory management, retailers' and customers' demands, stochastic nature of the SC, logistic issues, etc (Pinto, 2004; Cakravastia *et al.*, 2002; Giannoccaro and Pontrandolfo, 2002). The ultimate objective of the SCN optimization is to minimize the overall system wide cost while customer service is kept at pre-specified level (Truong and Azadivar, 2003) and maximizing profitability of not just one entity but rather all the entities in SCN, Pinto (2004). Moreover objective functions of the entities in the SCN are always conflicting in nature. For instance, the objective of marketing is to maximize customer service level and sales volume, which is in conflict with the objectives of production and distribution which is to produce and ship products of higher quality at minimum cost and in adequate amount. Also raw materials procurement decisions are aimed at minimizing cost of goods while production and distribution decisions are based on maximum output from plants with minimum production cost and demand. In such a scenario, there is a need to present a vector of options to the decision makers, so that the final decisions should be taken after taking a total balance over all criteria into account i.e. trade-offs (Latha Shankar *et al.*, 2013).

Multi-Objective Optimization (MOO) involves simultaneous optimization of problems with at least two objective functions which are conflicting in nature. In MOO, there is no single optimum solution, but a number of them exist that are optimal called Pareto fronts. There are several approaches of finding the Pareto fronts in MOO models. One can aggregate all the objectives into a single objective by scalarization using the weighted sum, distance functions, goal programming and $\varepsilon-$constraint (Konak *et al.*, 2006; Coello, 1999; Zitzler *et al.*, 2004. In contrast evolutionary algorithms are able to generate and maintain multiple solutions in one simulation. Commonly used evolutionary algorithms include genetic algorithms which have been shown to give a good approximation of the Pareto fronts (Kalyanmoy *et al.*, 2002; Konak *et al.*, 2006).

The supply chain network design problem has been optimized as single objective problem (Cohen and Lee, 1988; Arntzen *et al.*, 1995). However the results obtained from such models were not a true wholesomely optimal and to an extent misleading since an optimal solution for a specific scenario is not optimal to the other entities' objectives.

Several studies have been undertaken where the SCN problems have been optimized as MOO problems and the results are more than encouraging: Amodeo *et al.* (2007) optimized a supply chain as a multi-objective optimization using genetic algorithms and simulation model. Objectives considered were: minimizing inventory cost and maximizing service level. It was concluded that this approach obtained inventory policies better than the ones used in practice then with reduced costs and improved service level. Also Erol and Ferrell Jr (2004) proposed a model that assigned suppliers to warehouses and warehouses to customers and used multi-objective optimization modeling framework for minimizing system-wide costs and maximizing customer satisfaction. Chen and Lee (2004) developed a multi-product, multi-stage, and multi-period scheduling model for a multi-stage SCN with uncertain demands and prices. They simultaneously optimized the following objectives: fair profit distribution among all participants, safe inventory levels, maximum customer services, and robustness of decision to uncertain demands.

There has been a growing interest of using evolutionary algorithms to solve multi-objective optimization problems recently (Deb, 2001; Pinto, 2004; Farahani and Elahipanah, 2008). Different models have been developed with different objective functions where evolutionary algorithms have been used to find Pareto fronts. Sabri and Beamon (2000) developed an integrated multi-objective supply chain model for strategic and operational under certainties of products-delivery and demands. Similarly Melachrinoudis *et al.* (2005) worked on a bi-objective optimization with cost minimization and service level maximization as objectives. Pinto (2004) and Altiparmak *et al.* (2006) independently proposed a solution procedure based on genetic algorithms to find the Pareto optimal solutions for supply chain design problem.

Farahani and Elahipanah (2008) set up a bi-objective model for the distribution network of a supply chain which produces one product in a three echelon supply chain design. The objectives were: minimizing costs and minimizing backorders and surpluses of products. The Pareto optimal were found by using mixed integer programming by applying non-dominated sorting genetic algorithms.

Most of the existing models places much emphasis on the optimization of location and allocation decisions while ignoring important aspects such as capacities and technology aspects of the manufacturing facilities which also affects the production and distribution of products. In addition, many models have been implemented using genetic algorithms, genetic algorithms in particular. It should be noted that with large size problems they take long run time to find optimal solutions and sometimes they may converge towards a limited region of the Pareto fronts ignoring solutions that might be interesting. In such cases, there is need to guide the algorithm to converge towards desirable solutions if prior information is available.

This paper is organized as follows: section 2 gives the problem description in which we introduce the objective function and its assumptions followed by the mathematical models for the integrated supply chain. The proposed methodology is described in section 4 followed by computed results and discussion in section 5. Finally, the conclusion is given in section 6.

## 2. Problem Description

We consider a general supply chain network model consisting of suppliers, manufacturing plants, distribution centers (DC), and retailers as shown in Figure 1 below. The suppliers are companies from which raw materials are bought. We assume that each supplier supplies a specific raw material. In additional there are transporting vehicles of different capacities which transport the raw materials or products. The manufacturing plant is where the products are produced. The distribution centers are the warehouses of different capacities that store the products before they are delivered to retailers through use of vehicles. So we need to design a supply chain network that allows strategic decisions like selecting suppliers of raw materials quantities, determining the subsets of plants and distribution centers to be opened, and a distribution network strategy that satisfy all capacities and customer's demands in such a way that total cost and on time delivery delays are minimized.

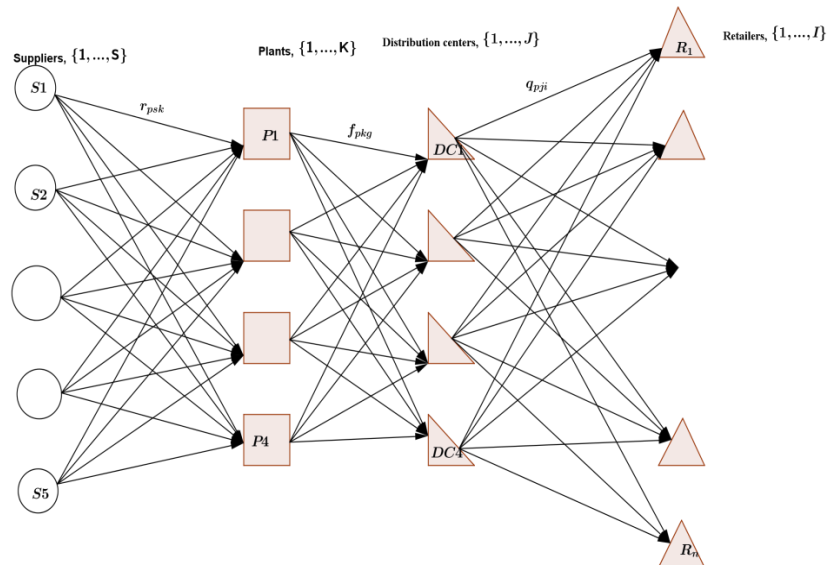

Figure 1. A Supply Chain Network Model

Ideally we would like to deliver the products at the right time the retailers have requested and in right amount as demanded. This helps to save storage cost and maximize customer services. So we require that every supplier must deliver the right amount of products at the right time and to the right place. Though products are stored prior to delivery, due to capacity constraints of the distribution centers and plants, it might not be possible to satisfy all the requests on time hence we either have to deliver earlier or late sometimes, incurring excess storage costs in the process. In this research we simultaneously optimize the following two objective functions: total cost of the supply chain which includes the cost of opening and running manufacturing plants and distribution centers, cost of buying and transporting raw materials, cost of transporting products from plants to DCs and from DCs to retailers, and cost of holding products at DCs. The other objective function is about on-time delivery (minimum delay) which involves earliness: representing the amount of products that are delivered prior to the due date, and tardiness: representing the amount of products that are delivered after their due date.

In terms of minimizing cost objective, it would be costly to supply products earlier to the distribution centers and deliver them late to retailers as we will incur holding costs. While for on time delivery objective we would want the warehouses to be supplied products earlier and have to store them up to an appropriate time so that we don't miss customer demands. As such we require a trade-off between delivering the products late or earlier and minimizing storage cost at the distribution center while maximizing customer service in terms of not missing customers' demands.

We considered the following assumptions:

- Number of retailers and suppliers and the capacity of suppliers are known.
- Number of plants and distribution centers and their maximum capacities are known.
- Demands of retailers are uncertain but can be determined from past history.
- We only consider that storage costs are incurred at distribution centers since that is where they may stay for long time before being delivered to retailers.
- Each retailer is served from one distribution center, while manufacturing plants delivers their products to many distribution centers.
- Each plant gets raw materials from all the suppliers.

### 3. The Mathematical Model

We use the following notations for our mathematical model formulation.

a) Notation for indices of the entities are as follows:
   $s$: supplier
   $k$: plant
   $j$: distribution center

$i$: retailer  
$p$: product type  
$t$: time  

b) Variables for quantities are as follows:

$r_{sk}$ : Amount of raw materials from supplier $s$ to plant $k$  
$f_{kj}$: Total number of products from plant $k$ to distribution center $j$  
$q_{ji}$: Total number of products from distribution center $j$ to retailer $i$  
$b_{pjt}$ : Total number of backorders of products at distribution center $j$ in time $t$ (tardiness)  
$z_{pjt}$ : Total number of products delivered to distribution center $j$ in time $t$ (earliness)  
$b_j$ : Binary variable to indicate whether a distribution center $j$ is open or closed  
$b_k$ : Binary variable to indicate whether a plant $k$ is open or closed  

c) The variable notation for model parameters are:

$D_k$: Capacity of plant $k$  
$SC_s$: Capacity of supplier $s$  
$d_i$: Quantity of demanded items from retailer $i$  
$d_{it}$: Quantity of demanded items from retailer $i$ on a given time $t$  
$DC_{tot}$: Total number of distribution centers  
$P_{tot}$: Total number of plants  
$v_j$: Fixed cost of operating distribution center $j$  
$g_k$: Fixed cost of operating plant $k$  
$c_{kj}$: Unit transportation of products from plant $k$ to distribution center $j$  
$t_{sk}$: Unit transportation cost of raw materials from supplier $s$ to plant $k$  
$c_s$ : Unit cost of raw material from supplier $s$  
$h_j$: unit holding cost at distribution center $j$  
$r_{ji}$: Unit transportation cost from distribution center $j$ to retailer $i$  
$Q_j$: Holding capacity of distribution center $j$  
$u$: utilization rate of raw materials per unit of each product  
$bl_{pjt}$: Maximum amount of permitted backorder for distribution center $j$ in time $t$.  

The mathematical formulation of the model is as follows:

minimize

$$f_1 = \sum_k g_k b_k + \sum_j v_j b_j + \sum_s \sum_k (c_s + t_{sk}) r_{sk} + \sum_k \sum_j c_{kj} f_{kj} + \sum_j \sum_t h_j b_{pjt} + \sum_j \sum_i r_{ji} q_{ji} \quad (0.1)$$

and minimize

$$f_2 = \sum_p \sum_j \sum_t (b_{pjt} + z_{pjt}) \quad (0.2)$$

Subject to the following constraints:

$b_{pjt} \leq Q_j$ : Warehouse's (distribution center) holding capacity.

$b_{pjt} \leq bl_{pjt}$ : Allowable backorder.

$\sum_i d_i \leq \sum_j Q_j$ : Capacity constraint for distribution centers.

$\sum_k f_{kj} \geq \sum_i q_{ji}$ : Amount of products sent to retailers is within the capacity of the distribution centers capacity.

$\sum_j q_{jt} = \sum_i d_{it}$ : Customer satisfaction.

$\sum_t d_{it} = d_i$ : Total demand of product by retailer $i$ over total time.

$\sum_k r_{sk} \leq SC_s$ : Supplier's capacity.

$u\sum_{j} f_{kj} \leq \sum_{s} r_{sk}$ : Supplier's capacity on production of products.

$u\sum_{j} f_{kj} \leq D_k$ : Plant's capacity for each product.

$\sum_{j=1}^{DC_{tot}} b_{ji} = 1, i = 1, 2, ..., N$ : Each customer is served by a single DC.

$b_j \in \{0,1\}$, $b_k \in \{0,1\}$ : Binary variables

$f_{kj}, r_{sk}, g_k, v_j, t_{sk}, c_{kj}, bl_{pjt}, Q_j, b_{pjt}, z_{pjt}, SC_s, D_k, d_{it} \geq 0$

## 4. Multi-Objective Evolutionary Algorithms (MOEA)

In general a MOO problem has the form

$$\min \boldsymbol{F} = \{f_1(\boldsymbol{x}), f_2(\boldsymbol{x}), ..., f_N(\boldsymbol{x})\}$$
$$\text{subject to } g_i(\boldsymbol{x}) \leq 0, \ i=1,...,M$$
$$h_j(\boldsymbol{x})=0, \ j=1,...K$$
$$\boldsymbol{x}_l \leq \boldsymbol{x} \leq \boldsymbol{x}_u$$

In MOO, mathematically, a feasible solution $x_1$ is said to dominate another feasible solution $x_2$ if and only if $f_i(x_1) \leq f_i(x_2)$, $i = 1, 2, ..., N$ and $f_j(x_1) < f_j(x_2)$ for at least one objective function $j$. A solution that is not dominated by any other solution is called Pareto-optimal. The set of all non-dominated solutions is called Pareto-Optimal set and the corresponding objective values for the Pareto-optimal set in the objective space is called Pareto-fronts. Hence the ultimate goal in MOO is to find the Pareto fronts which give trade-offs among the objectives being optimized. Pareto fronts cannot be improved with respect to any objective without worsening at least one other objective function.

Identifying the entire Pareto-optimal set, for many multi-objective problems is practically impossible due to its size, in addition, for combinatorial problems like supply chain network, proof of optimality is computationally infeasible (Konak *et al.*, 2006). Therefore we investigate a set of solutions that approximate the Pareto-optimal set as much as possible. Over the past decades, several multi-objective evolutionary algorithms (MOEA) have been proposed by different researchers to solve MOO problems (Deb *et al.*, 2002). The advantages of MOEAs is that they deal with a set of possible solutions called populations which enables finding of the entire set of Pareto fronts and maintaining them in one run of the algorithm which is one of the necessary element in multi-objective optimization (Deb, 2001; Donoso and Fabregat, 2010) . In addition they are less susceptible to the shape or continuity of the Pareto fronts (Coello, 1999; Deb *et al.*, 2002). Evolutionary algorithms originate in Darwin's theory of evolution and were first developed as solution for optimization problems by Holland (1975). Simulating the biological evolution process, the algorithms use structure or individuals to solve the problem (Donoso and Fabregat, 2010; Konak *et al.*, 2006). The algorithm consists of individuals called chromosome, the population of the individuals, the fitness of each of the individuals that is its phenotype, and the genetic operators namely mutation, selection and cross over. Each chromosome is made up of discrete units called genes that controls one or more features of the chromosome, and corresponds to a unique solution. In crossover, generally two chromosomes referred as parents combine to form new chromosomes called offspring while mutation introduces random change into characteristics of chromosomes which helps to re-introduce genetic diversity back into the population so as to avoid converging to a local optimal. Selection is the process of choosing chromosomes among the population to form the next generation. Some methods of selection include proportional selection, ranking and tournament selection and $(\mu + \lambda)$ selection. Evolutionary algorithms run in polynomial time, so even though there is no guarantee that an optimal value will necessarily be found, we can find a good approximate value that we can actually use (Donoso and Fabregat, 2010; Deb, 2001).

The first MOEA was vector evaluated genetic algorithm (VEGA) proposed by Schaffer (1984), which was an extension of simple genetic algorithm to accommodate vector-valued fitness measures followed by multi-objective genetic algorithm (MOGA) proposed by Fonseca and Fleming (1993). Both VEGA and MOGA had no mechanisms to preserve best fit individuals that should have been passed to the next generation which enhances convergence to optimal solutions. Since then, several MOEAs have been developed namely niched pareto genetic algorithm, random weighted genetic algorithm (RWGA), non-dominated sorting algorithm (NSGA), strength Pareto evolutionary algorithm (SPEA), and fast non dominated sorting algorithm (NSGA-II) (Konak *et al.*, 2006; Deb *et al.*, 2002; Srinivas and Deb, 1994; Coello, 1999). They all incorporate crossover, selection, mutation and replacement, but they differ in the way they implement them and how they deal with multiple objectives and preservation of diversity in the solutions.

In this research we make use of NSGA-II genetic algorithm proposed by Deb *et al.* (2002) which was aimed at overcoming problems associated with NSGA. NSGA developed by Srinivas and Deb (1994) was based exactly on MOGA except for the specification of the sharing parameter. However it has higher computational complexity of $O(MN^3)$ for $M$ objectives and a population size of $N$ lacked elitism and required specification of a sharing parameter for the fitness value. Unlike NSGA, NSGA-II implements elitism which is a mechanism that allows preservation of best fit individuals in a population to ensure good fitness already obtained does not get lost in subsequent generations. It has a reduced computational complexity of $O(MN^2)$. NSGA-II uses tournament selection whereby a group of individuals takes part in a tournament and the winner is judged by fitness levels that each individual brings to the tournament. It has been shown that it outperforms most contemporary MOEAs like SPEA, PAES, etc and its performance has been tested on several test problems where it has given accurate results in generating Pareto fronts (Deb *et al.*, 2002; Pinto, 2004).

Our choice of NSGA-II for the research is based on the fact that it is an elitist and fast strategy, modular and flexible, emphases on the non dominated solutions, can be applied to a large wide of problems and the availability of global optimization toolbox in MATLAB that implements it so that we do not re-invent the wheel.

The main principle in NSGA-II is that we categorize all solutions into ranks and for each individual solution we calculate the crowding distance. To determine the ranking, two entities are calculated for each solution $p$ namely $n_p$ domination count representing the number of solutions that dominate solution $p$, and $s_p$ set of all solutions $q$ dominated by solution $p$. All solutions $p$ with $n_p = 0$ are identified and ranked as 1, the first non-domination front. Thereafter, we visit each member $q$ of $s_p$ where $n_p = 0$ and reduce its domination count by 1. If $n_q = 0$ then solution $q$ is ranked as 2, the second non-domination front. We continue in this way until all ranks are identified. Solutions with rank 1 are non-dominated and represent the best. Then we calculate the overcrowding distance of each solution which is the measure of the diversity of solutions along the front of the non dominated solutions and enables the algorithm to obtain uniformly distributed solutions over the true Pareto fronts. The crowding distance of each solution $I[i]_{distance}$ is calculated as follows:

Let $n=|I|$, $I$ is set of solutions. For each $i \in I$, assign $[I]_{distance} = 0$.

For each objective $m$ sort the set $I$ in ascending order based on the objective values.

For each objective $m$, assign large distances to the boundary solutions, $I[1]_{distance} = I[n]_{distance} = \infty$. For all solutions $i=2$ to $n-1$, then

$$I[i]_{distance} = I[i]_{distance} + \frac{(I[i+1]_m - I[i-1]_m)}{f_m^{\max} - f_m^{\min}}, \text{ where}$$

$I[i]_m$ : $m^{th}$ objective function value of the $i^{th}$ individual in the set $I$,

$f_m^{\max}$ : maximum value of the $m^{th}$ objective function,

$f_m^{\min}$ : minimum value of the $m^{th}$ objective function.

NSGA-II uses a two step selection process, which combines binary tournament and $(\mu + \lambda)$ selection (Reed *et al.*, 2003). The binary tournament allows only the fittest individuals to be placed in the mating population while $(\mu + \lambda)$ selection, selects which of the parents or children should form the next generation. To enable elitism, individuals for the next generation are selected from both the parents and children based on their non-domination front and crowding distance in case of tier until $N$ individuals. Here are the steps involved in the algorithm:

**Step1**: Generate a parent population $P_t$ randomly and create an offspring population $Q_t$ from $P_t$ using genetic operators: selection, crossover and mutation each of size $N$.

**Step 2**: Combine $P_t$ and $Q_t$ into a single population set $S_t$ where $|S_t| = 2N$.

**Step3**: Classify and rank each individual in $S_t$.

**Step 4:** Select the best $N$ individuals from $S_t$ based on ranks and crowding distance to form a parent population $P_{t+1}$ for the next generation.

**Step 5:** Repeat steps 1-4 until the termination criteria is met.

For each of plant, distribution center and retailer, we dedicate a section in the chromosome for its representation as in Figure 2.

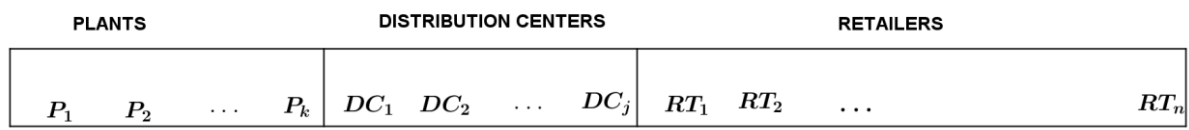

Figure 2. An example chromosome

The variables $P_k$ represents amount of raw materials received on each plant $P$, while $DC_j$ are the amount of products that are produced and transported to distribution center $DC$. $RT_n$ is amount of products received on each retailer $RT$.

## 5. Results and Discussion

In this research, data was collected from SBC Tanzania which is a supply chain company that produces and sells soft drinks to both retailers and wholesalers (SBC Tanzania, 2010). The simulation was run in MATLAB 2013A with the following parameters: population size is 1290 (default) which is initially uniformly created, selection is by tournament with crossover probability of 0.6. Mutation is at 0.01. The mutation and crossover probabilities are chosen to tally with the recommendation of Reed *et al.* (2003) on choices of parameters for NSGA-II after we had run the simulation for several values. The Pareto fronts are shown in the figure below for a maximum generation of 500 and the default maximum value in MATLAB :

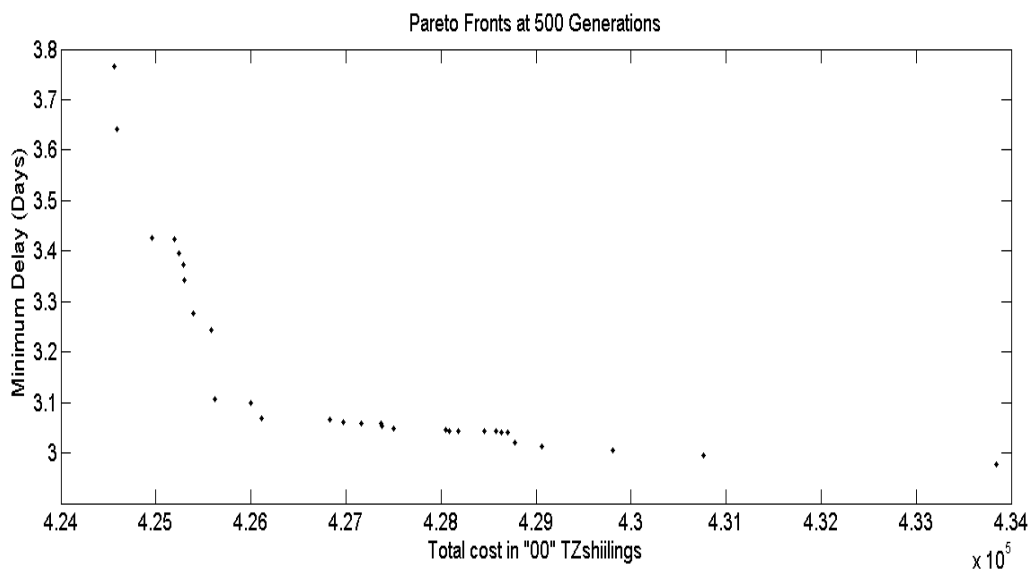

Figure 3. Pareto Fronts at 500 Generations

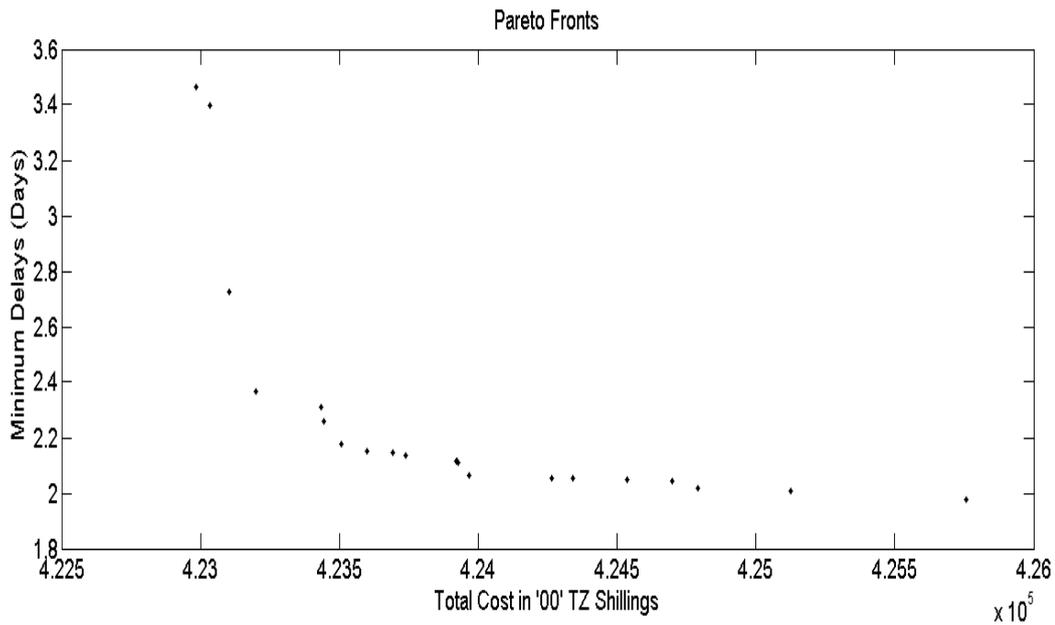

Figure 4. Pareto Fronts Using Maximum Generations

The objective function values are shown in the Table 1 and Table 2 below:

Table 1. Pareto Fronts at 500 Generations

| Total Cost (TZ shillings/week) | Minimum Delays (Days) | Total Cost (TZ shillings/week) | Minimum Delays(Days) |
|---|---|---|---|
| 42496345 | 3.43 | 42697666 | 3.06 |
| 42529667 | 3.34 | 42750136 | 3.05 |
| 43384315 | 2.98 | 42716485 | 3.06 |
| 42524359 | 3.40 | 42599589 | 3.10 |
| 43075927 | 2.99 | 42906090 | 3.01 |
| 42804906 | 3.05 | 42683399 | 3.07 |
| 42529667 | 3.34 | 42737024 | 3.06 |
| 42558216 | 3.24 | 42863937 | 3.04 |
| 42529430 | 3.37 | 42539146 | 3.28 |
| 42610867 | 3.07 | 42808494 | 3.04 |
| 42519600 | 3.42 | 42980958 | 3.00 |
| 42818656 | 3.04 | 42459411 | 3.64 |
| 42846044 | 3.04 | 42858301 | 3.04 |
| 42529667 | 3.34 | 42737934 | 3.05 |
| 42529667 | 3.34 | 43384315 | 2.98 |
| 42562288 | 3.11 | 42456237 | 3.77 |
| 42870315 | 3.04 | 42558216 | 3.24 |
| 42877834 | 3.02 | | |

Table 2. Pareto Fronts at Maximum Generations

| Total Cost (TZ Shillings/week) | Minimum Delay (Days) | Total Cost (TZ Shillings/week) | Minimum Delay (Days) |
| --- | --- | --- | --- |
| 42392602 | 2.11 | 42343410 | 2.31 |
| 42512544 | 2.01 | 42453549 | 2.05 |
| 42479149 | 2.02 | 42373738 | 2.13 |
| 42426517 | 2.06 | 42359834 | 2.15 |
| 42319822 | 2.37 | 42433977 | 2.05 |
| 42344611 | 2.26 | 42469855 | 2.04 |
| 42396879 | 2.06 | 42350870 | 2.18 |
| 42369365 | 2.14 | 42303267 | 3.40 |
| 42392100 | 2.12 | 42310332 | 2.73 |
| 42575749 | 1.98 | 42392297 | 2.11 |
| 42298487 | 3.46 | | |

With this vector of optimal choices offering several trade-offs, the decision maker can thus intelligently decide based on expertise to evaluate those solutions that would be most beneficial to the company from a system's optimization perspective.

On the results itself, it was found that there is an inverse relationship between the total cost and the allowable delay in both simulations though with a slight difference. If products are delayed for many days the total cost decreases and if demands are met with minimum delays the cost becomes bigger. The results can be interpreted as follows: to satisfy customers demands we need to have minimum delays in delivering goods once they have requested. Since goods are manufactured from plants and temporally stored in distribution centers then we incur several costs in terms of production as there is need to manufacture more goods so that they are readily available once they are ordered. In the same way these products may be kept for long time at distribution centers hence incurring more storage costs. There might also be a need for more transporting vehicles hence increasing transportation cost. With maximum delays these additional costs may be small or not there at all as we can deliver goods as long as they are available with no pressure.

From our perspective, we would advise the industry not to delay the deliveries with more than two days as it makes a good business sense and the earlier the delivery the better. And it is advisable most times to choose extremes when a decision cannot be made on which option to select.

The results obtained in this optimization ranges from 40-50 million TZ Shillings which is less compared to the actual amount that is spend by the company in a week which is 45-65 million TZ shillings on the same processes considered.

## 6. Conclusion

In this paper, we have developed a mathematical model for a distribution network in a three echelon supply chain that minimizes the system-wide costs and delays on delivery of products. The mathematical model is designed as a multi-objective optimization problem taking into account the two conflicting objective functions. We have found several options in which, for some cases we might opt for high costs to maintain customer satisfaction, whereas in other cases depending on the situation we might opt to save costs. We used NSGA-II through the Global Optimization Toolbox in MATLAB. The algorithm clearly provides the Pareto fronts which are efficient solutions so that decision makers can use in planning and designing to improve the supply chain services. This work demonstrates that evolutionary algorithms provide a successful way of dealing with multi-objective optimization problems and have potential to solve combinatorial problems. Though the MOEAs have proved to be efficient in solving multi-objective problems, they become unsuccessful as the number of objectives increases from two in the problem. The main principle of these algorithms is the emphasis on non-dominated solutions in a population; however, as the number of objectives increases most population members in a randomly created population tend to become non-dominated to each other resulting in the algorithm failing to introduce new population members in the new generation thereby causing a stagnant performance. For instance 10% and 90% of members in a population of size 200 are non-dominated in a 3 and 10 objective functions problem respectively (Deb, 2011; Reed *et al.*, 2003).

Future research should consider designing of a model where the company hires transporting vehicles rather than having their own as there are other costs like maintenance and wear and tear. In addition, it might be interesting to if reducing product's cost and allowing retailers to organize their own transport to collect products from distribution centers is optimal or not as compared to the current scenario where the company delivers products to retailers as the number of customers keep increasing.

**Acknowledgement**

Nelson is indebted to the Nelson Mandela African Institution of Science and Technology for financial and material support.